\documentclass[reqno]{amsart}
\usepackage{amssymb}
\usepackage{amsmath}
\usepackage{amsfonts}
\usepackage{graphicx}
\newtheorem{thm}{Theorem}
\newtheorem{lem}[thm]{Lemma}

\theoremstyle{definition}

\theoremstyle{definition}

\theoremstyle{definition}

\begin{document}
\title{Continuity of generalized wave maps on the sphere}

\author{Daniel {d}a Silva}
\address{Department of Mathematics, University of Rochester, Rochester, NY 14627}
\date{\today}

\begin{abstract}
We consider a generalization of wave maps based on the Adkins-Nappi model of nuclear physics.  In particular, we show that solutions to this equation remain continuous at the origin, which is a first step towards establishing a regularity theory for this equation.
\end{abstract}

\maketitle

\section{Introduction}
\label{intro}

Let $\phi: \mathbb{R}^{1+n} \to N$, where $(\mathbb{R}^{1+n}, g)$ is the $1 + n$ dimensional Minkowski spacetime with metric $\eta = \text{diag}(-1,1,\ldots,1)$, and $(N, h)$ is a Riemannian manifold.  We say that $\phi$ is a \emph{wave map} if it is a formal critical point of the action
\begin{equation}
\label{WML}
S = \frac{1}{2} \int f^{\mu \nu} \partial_{\mu} \phi^{i} \partial_{\nu} \phi^{j} h_{i j}(u) \ dx dt = \frac{1}{2} \int f^{\mu \nu} S_{\mu \nu}\ dx dt.
\end{equation}
Here and throughout, we use the Einstein summation convention.  Critical points of \eqref{WML} satisfy the wave maps equation,
\begin{equation}
\label{WM}
\Box \phi^i + \Gamma^{i}_{j k}(\phi)\partial^{\alpha}\phi^{j} \partial_{\alpha} \phi^{k} = 0.
\end{equation}
Here the $\Gamma^{i}_{j k}$ are the Christoffel symbols corresponding to $h$.

Of particular interest to physicists is the case where $n = 3$ and $N = \mathbb{S}^{3}$.  This is known to particle physicists as the \emph{nonlinear $\sigma$ model}, first proposed by Gell-Mann and L\'evy as a model for interactions of particles known as \emph{pions}.  This model, however, proved to be inadequate, due to a result known as Derrick's theorem \cite{Derrick1964}, which implies the non-existence of topological solitons (static, localized solutions) for the nonlinear $\sigma$ model.  As a result, several alternative models were proposed, each with a different approach to avoiding Derrick's theorem.

One of these models is known as the Adkins-Nappi model \cite{Adkins1984}.  In this model, we again assume $n=3$ and $N = \mathbb{S}^3$.  Derrick's theorem is avoided, however, by introducing additional terms to the action, which correspond to interactions between pions and other particles known as \emph{vector mesons}.  The new action is given by
\begin{equation}
\label{ANL}
S = \frac{1}{2} \int f^{\mu \nu} S_{\mu \nu} \ dxdt + \frac{1}{4} \int F^{\mu\nu}F_{\mu\nu}\ dg - \int A_{\mu}j^{\mu}\ dxdt.
\end{equation}
Here, $A = A_{\mu} dx^{\mu}$ is a gauge potential representing the vector mesons, $F$ is its associated electromagnetic field, given by
\[
F_{\mu\nu} = \partial_{\mu} A_{\nu} - \partial_{\nu} A_{\mu},
\]
while $j$ is the baryonic current
\[
j^{\mu} = \frac{1}{24\pi^2} \epsilon^{\mu\nu\rho\sigma}\partial_{\nu}\phi^i \partial_{\rho} \phi^j \partial_{\sigma} \phi^k \epsilon_{ijk}.
\]

This model was recently studied by Geba and Rajeev, who began a program to prove the existence of global smooth solutions.  As a first step, they focused on the equivariant case, for which solutions take the special form
\[
\phi (t,x) = \phi (t,r,\omega) = (u(t,r), \omega),
\]
where $\omega \in \mathbb{S}^2$.  A simple computation will show that the Euler-Lagrange equations for $\phi$ and $A$ decouple, and that the equations for the components of $\omega$ are trivially satisfied.  The equation for $u$ becomes
\begin{equation}
\label{AN}
u_{tt} - u_{rr} - \frac{2}{r}u_r + \frac{\sin 2u}{r^2} + \frac{(u - \sin u \cos u)(1- \cos 2u)}{r^4} = 0.
\end{equation}
Geba and Rajeev showed in \cite{Geba2010} that local solutions to the Cauchy problem for \eqref{AN} are continuous under the assumption of smooth data of finite energy.  The goal of the present work is to generalize this result for more general nonlinearities, and for all dimensions $n \geq 2$.  

To motivate our result, consider \eqref{WM} for $N = \mathbb{S}^n$, with $n \geq 2$.  If we also make the ansatz $\phi (t,r, \omega) = \big(u(t,r), \omega\big)$, with $\omega \in \mathbb{S}^n$, it is easily seen that, for the $u$ coordinate, \eqref{WM} becomes
\[
u_{tt} - u_{rr} - \frac{n-1}{r}\, u_r + \frac{n-1}{2}\, \frac{\sin 2u}{r^2} = 0.
\]
We will consider a generalization based of this equation based on \eqref{AN}.  In the case $n=3$, the additional terms of the Adkins-Nappi Lagrangian introduced nonlinearities of the form
\[
\frac{f(u) f^{\prime}(u)}{r^{\alpha}},
\]
where $f(0) = 0$.  This leads us to consider the Cauchy problem
\begin{equation}
\label{GWM}
\begin{cases}
u_{tt} - u_{rr} - \dfrac{n-1}{r}\, u_r + \dfrac{n-1}{2}\, \dfrac{\sin 2u}{r^2} + \dfrac{f(u)f^{\prime}(u)}{r^{\alpha}} = 0 \\
u(t_0) = u_0 \\
u_t(t_0) = u_1,
\end{cases}
\end{equation}
where $u_0$, $u_1$ are $C^{\infty}$ functions.

\section{Formulation and Statement of Results}

We begin our analysis with the observation that equation \eqref{GWM} is a semilinear hyperbolic equation.  As such, the local existence theory assures us that smooth solutions exist, at least for a finite time $T_0$.  Our ultimate goal is to show that the conserved energy for \eqref{GWM}, given by
\[
E(u) = \frac{1}{2} \int_{\mathbb{R}^n} u_t^2 + u_r^2 + (n-1) \frac{\sin^2 u}{r^2} + \frac{f^2(u)}{r^{\alpha}}\ dx,
\]
does not concentrate.  Due to the symmetry of the problem, if the energy does concentrate, then it can only do so at $r = 0$.  This indicates that the origin requires special attention.

As a first step towards proving non-concentration, we will show that solutions remain continuous at the origin.  Since \eqref{GWM} is invariant under translations and reflections in time, we can consider the Cauchy problem where the data is given at time $t_0 = T_0 > 0$, and the possible blow-up occurs when $t = 0$.

To understand how such a result will be proven, we must examine the conserved energy for \eqref{GWM} more closely.  Due to the presence of the terms
\[
\int_{\mathbb{R}^n} \frac{\sin^2 u}{r^2}\ dx \quad
\textrm{and} \quad
\int_{\mathbb{R}^{n}} \frac{f^2(u)}{r^{\alpha}}\ dx
\]
in the expression for the energy, we must impose the boundary condition $u(t,0) = 0$ for all $t$ to ensure that we have solutions of finite energy.  Based on this, we see that we need only show that
\[
\lim_{(t,r) \to (0,0)} u(t,r) = 0.
\]
Due to the finite speed of propagation, it suffices to consider $(t,r)$ in the forward light cone
\[
\Omega = \{ (t,r):\ 0 \leq r \leq t,\ 0 < t \leq T\}.
\]
Thus, in the subsequent analysis, when we say that $u$ is a smooth local solution to \eqref{GWM}, we mean that $u$ is a solution to \eqref{GWM} which is smooth in $\Omega$.

With these facts in mind, we can now state the main theorem that will be proved.
\begin{thm}
\label{bigthm1}
If $u$ is a smooth local solution of finite energy to \eqref{GWM}, where  
\[\alpha \geq \max\{ 2(n-1), n+1 \}
\]  
and
\[
f(0) = 0, \quad f^{\prime}(0) \neq 0,\quad f(u) f^{\prime}(u)  u \geq 0, \quad\int_0^\infty |f(v)|\ dv = \infty,
\]
then
\begin{equation}
\label{energyan}
\lim_{T\to 0_+} \int_{\Sigma_T} \left(1 - \frac{r}{t}\right)\left(u_t -u_r\right)^2 + (u_t + u_r)^2 + \frac{\sin^2 u}{r^2} + \frac{f^2(u)}{r^\alpha} = 0.
\end{equation}
Moreover, $u$ is continuous at the origin:
\[
\lim_{(t,r)\to (0,0)} u(t,r)=0.
\]
\end{thm}

\section{Proof of Theorem \ref{bigthm1}}

The proof of Theorem \ref{bigthm1} is based on the formalism of Shatah and Tahvildar-Zadeh in \cite{Shatah1992} for the study of wave maps.  This formalism makes use of the method of multipliers, first introduced by Friedrichs and later developed by Morawetz in \cite{Morawetz1975}.  We will begin by introducing notation and basic results, which will be used in subsequent sections to prove the main portions of the result.

\subsection{Preliminaries}
\begin{figure}[h]
\centering 
\includegraphics[scale=.3]{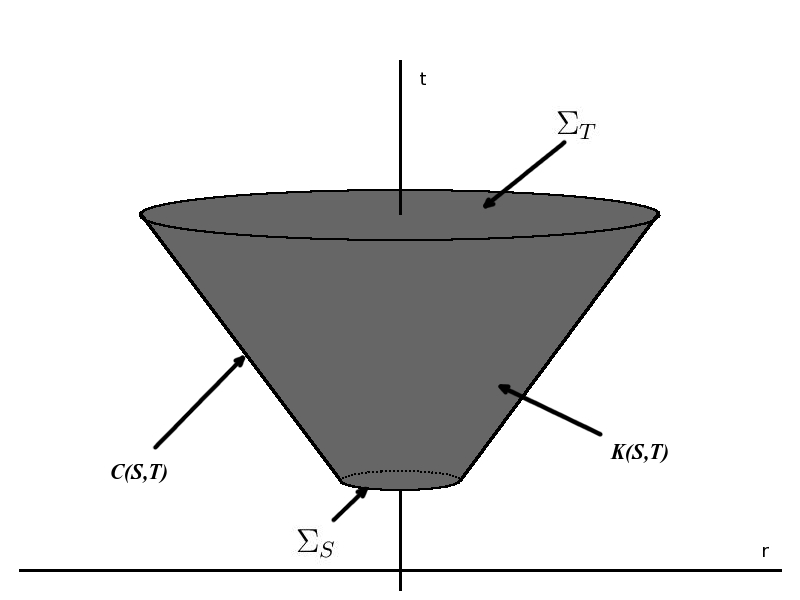} 
\caption{The truncated forward light cone.} 
\label{fig1} 
\end{figure}

We start by introducing the notation
\[
\begin{aligned}
e^{\pm}(t,r) & = \frac{1}{2}\left(u_t^2(t,r) + u_r^2(t,r)\right) \pm  \frac{n-1}{2}\, \frac{\sin^2 u(t,r)}{r^2} \pm \frac{f^2(u(t,r))}{2r^\alpha}, \\
m(t,r) & = u_t(t,r)u_r(t,r). 
\end{aligned}
\]
The quantity $e^+$ is often referred to as the energy density, while $m$ is often called the momentum density.  Multiplying \eqref{GWM} by $r^{n-1} (a u_t + b u_r + c u)$ and rewriting the result as a divergence, we obtain the main differential identity which will be used in our computations.
\begin{lem}
If $u$ is a smooth solution to \eqref{GWM}, then it satisfies the following identity:
\begin{equation}
\label{gendiffid}
\begin{aligned}
& \partial_t\left[r^{n-1}(a \, e^+ + b \, m + c \, u \, u_t) \right] - \partial_r\left[r^{n-1}(a \, m + b \, e^- + c \, u \, u_r) \right] \\
& \qquad \qquad = \frac{r^{n-1}}{2}\left( a_t - b_r - (n-1)\frac{b}{r} + 2c \right) u_t^2 \\
& \qquad \qquad \ + \frac{r^{n-1}}{2}\left( a_t - b_r + (n-1)\frac{b}{r} - 2c \right) u_r^2 \\
& \qquad \qquad \ + r^{n-1} \left( a_t + b_r + (n-3)\frac{b}{r} \right)\frac{n-1}{2}\,  \frac{\sin^2 u}{r^2} \\
& \qquad \qquad \ + r^{n-1} \left( a_t + b_r + \left( n - 1- \alpha\right)\frac{b}{r}\right) \frac{f^2(u)}{2r^{\alpha}} \\
& \qquad \qquad \ + r^{n-1}(b_t - a_r) u_t u_r + r^{n-1} u  (c_t u_t - c_r u_r) \\
& \qquad \qquad \ - r^{n-1}c u\left(\frac{n-1}{2}\, \frac{\sin 2u}{r^2} + \frac{f^{\prime}(u)f(u)}{r^{\alpha}}\right).
\end{aligned}
\end{equation}
\end{lem}
This identity will be used with integrals over the following subsets of $\mathbb{R}^{1+n}$ (see Figure \ref{fig1}).
\[
\begin{aligned}
K(S,T) & = \{ (t,r):\ 0 \leq r \leq t,\ 0 < S \leq t \leq T \leq T_0\} \\
\Sigma_T & = \{ (T,r):\ 0 \leq r \leq T \} \\
C(S,T) & = \{ (t,t):\ 0 < S \leq t \leq T \leq T_0\}
\end{aligned}
\]

Next, we would like to prove some preliminary lemmas which will be used below.  Consider $(a,b,c)=(1,0,0)$ in \eqref{gendiffid}, from which we deduce the energy differential identity
\begin{equation}
\label{energy}
\partial_t\left(r^{n-1} e^+ \right) - \partial_r\left(r^{n-1} m \right)=0.
\end{equation}
This will be used to relate the local energy
\[
E(T) = \int_{\Sigma_{T}} e^{+}(t,r),
\]
and the flux,
\[
F(S,T) = \frac{1}{\sqrt{2}} \int_{C(S,T)} (e^{+}+ m)(t,r),
\]
via the following lemma.
\begin{lem}
\label{mono}
For smooth local solutions of finite energy and $S \leq T \leq T_0$, we have 
\begin{equation}
E(T) - E(S)=F(S,T),
\label{flux}
\end{equation}
which implies that the energy is monotone, while 
\begin{equation}
\lim_{T \to 0^{+}} F(T) = 0,
\label{decayflux}
\end{equation}
with 
\[
F(T) = \lim_{S \to 0^+} F(S,T).
\]
\end{lem}
\begin{proof}
We can obtain \eqref{flux} immediately by integrating the energy differential identity on the frustum $K(S,T)$, then applying the divergence theorem.  For the second statement, first we note that
\[
e^{+}+ m= \frac{(u_t+u_r)^2}{2} + \frac{n-1}{2}\, \frac{\sin^2 u}{r^2} + \frac{f^2(u)}{2r^\alpha} \geq 0.
\]
It follows that $F(S,T) \geq 0$ for all $T \geq S \geq 0$.  Combining this fact with \eqref{flux}, we may conclude that 
\[
E(T) \geq E(S).
\]

For the final statement, we note that the monotonicity of $E$, combined with \eqref{flux}, tells us that $F(S,T)$ is uniformly bounded in $S$.  Since the integrand in the definition of $F(S,T)$ is positive, we can conclude that the limit
\[
\lim_{S \to 0^+} F(S,T)
\]
exists.  By dominated convergence, it follows that
\[
\lim_{T \to 0^+} F(T) = 0.
\]
\end{proof}


\subsection{Boundedness and Continuity Arguments}
We start this section by introducing the Bogomolny functional
\[
I(w) = \int_0^w |f(v)|\ dv,
\]
which, based on the hypotheses on $f$, listed in Theorem \ref{bigthm1}, has the following properties
\begin{itemize}
\item $I=I(w)$ is continuous;
\item $I(w)\cdot w  \geq 0$, with equality only at $w=0$;
\item $|I(w)| \rightarrow \infty$ as $|w| \rightarrow \infty$.
\end{itemize}

For $u$ a local smooth solution of \eqref{GWM}, the fundamental theorem of calculus leads us to
\[
\begin{aligned}
I(u(t,r)) & = \int_0^r I^{\prime}(u(t,s)) u_r(t,s)\ ds \\
& = \int_0^r |f(u(t,s))| u_r(t,s)\ ds,
\end{aligned}
\]
due to $u(t,0)=0$. Following this with an application of the Cauchy-Schwarz inequality and relying on $\alpha\geq 2(n-1)$, we deduce
\begin{equation}
\label{bogo}
\begin{aligned}
\left| I(u(t,r)) \right| & = \left| \int_0^r |f(u(t,s))| u_r(t,s)\ ds \right| \\
& \leq \left(\int_0^r \frac{f^2(u(t,s))}{s^\alpha}s^{\alpha-2(n-1)}\ s^{n-1}ds \right)^{\frac{1}{2}} \left(\int_0^r u_r^2(t,s)\ s^{n-1}ds \right)^{\frac{1}{2}} \\
& \leq r^{\frac{\alpha-2(n-1)}{2}} \left(\int_0^r \frac{f^2(u(t,s))}{s^\alpha}\ s^{n-1}ds \right)^{\frac{1}{2}} \left(\int_0^r u_r^2(t,s)\ s^{n-1}ds \right)^{\frac{1}{2}} \\
& \lesssim r^{\frac{\alpha-2(n-1)}{2}} E^{\frac{1}{2}}(t) \left(\int_0^r \frac{f^2(u(t,s))}{s^\alpha}\ s^{n-1}ds \right)^{\frac{1}{2}}  \lesssim r^{\frac{\alpha-2(n-1)}{2}} E(t).
\end{aligned}
\end{equation}

This inequality will serve several purposes below.  First, we will use it to show that $u$ is bounded in the forward light cone.
\begin{lem}
If $u$ is a local smooth solution to \eqref{GWM} of finite energy, then \\ $u\in L^{\infty}(K(0,T_0)).$
\end{lem}
\begin{proof}
Taking advantage of the monotonicity of the energy proved in Lemma \ref{mono}, we can write for all $(t,r)\in K(0,T_0)$:
\[
\left| I(u(t,r)) \right| \lesssim r^{\frac{\alpha-2(n-1)}{2}} E(t) \leq  T_0^{\frac{\alpha-2(n-1)}{2}} E(T_0) .
\] 
However, by construction, we know that $|I(w)| \rightarrow \infty$ as $|w| \rightarrow \infty$. This obviously implies the desired conclusion.
\end{proof}

Secondly, if we write \eqref{bogo} for the case $n=2$, we obtain
\[
\left| I(u(t,r)) \right| \lesssim r^{\frac{\alpha-2}{2}} E(t),
\] 
which proves the continuity of the solution, based on $\alpha \geq \max\{ 2(n-1), n+1 \}=3$ and the properties of the $I$ functional.  The case $n \geq 3$, however, is more subtle.  For this case, we will require that part of the energy does not concentrate, as indicated in the following lemma.
\begin{lem}
\label{biglem}
If $u$ a local smooth solution to \eqref{GWM} of finite energy and
\begin{equation}
\label{zoom}
\lim_{T \to 0^+} \int_{\Sigma_T} \frac{f^2(u)}{r^\alpha}=0,
\end{equation}
then
\[
\lim_{\substack{(t,r) \to (0,0) \\ (t,r) \in K(0, T_0)}} u(t,r)=0.
\]
\end{lem}
\begin{proof}
Using the last line of \eqref{bogo} and the monotonicity of the energy, we infer
\[
\aligned
\left| I(u(t,r)) \right| &\lesssim r^{\frac{\alpha-2(n-1)}{2}} E^{\frac{1}{2}}(t) \left(\int_0^r \frac{f^2(u(t,s))}{s^\alpha}\ s^{n-1}ds \right)^{\frac{1}{2}}\\  &\lesssim r^{\frac{\alpha-2(n-1)}{2}} E^{\frac{1}{2}}(T_0) \left(\int_{\Sigma_t} \frac{f^2(u)}{s^\alpha} \right)^{\frac{1}{2}},
\endaligned
\]
which, based on the hypothesis on $\alpha$ and properties of the Bogomolny functional, finishes the proof.
\end{proof}

Thus, for $n \geq 3$, Theorem \ref{bigthm1} will follow from Lemma \ref{biglem}.  The remainder of this paper will be devoted to proving the limit in \eqref{zoom}.

\subsection{Main Argument}
We will now complete the proof of Theorem \ref{bigthm1} by proving the limit in \eqref{energyan}.  We start by writing \eqref{gendiffid} for $a=t$, $b=r$, and $c=(n-1)/2$:
\begin{equation}
\begin{aligned}
\label{3diffid}
&\partial_t\left[r^{n-1}\left(t e^+ + rm + \frac{n-1}{2}u u_t\right) \right] - \partial_r\left[r^{n-1}\left(t m + r e^- + \frac{n-1}{2} u u_r\right) \right]  \\
& \qquad \qquad  = r^{n-1}\left[\frac{(n-1)^2}{2}\,\frac{\sin^2 u}{r^2} + (n+1- \alpha)\frac{f^2(u)}{2r^{\alpha}}\right. \\ 
& \left. \qquad \qquad \quad - \frac{n-1}{2} u \left(\frac{n-1}{2}\, \frac{\sin (2u)}{r^2} + \frac{f^{\prime}(u)f(u)}{r^{\alpha}}\right)\right]
\end{aligned}
\end{equation}
Integrating this over the region $K(S,T)$ yields
\begin{equation}
\begin{aligned}
\label{energyint}
& \int_{K(S,T)} \frac{(n-1)^2}{4}\,\frac{u f^{\prime}(u)f(u)}{r^{\alpha}} + (\alpha - (n+1))\frac{f^2(u)}{2r^{\alpha}} \\
& + \int_{\Sigma_T} Te^+ + r u_t u_r + \frac{n-1}{2} u u_t = \int_{\Sigma_S} Se^+ + r u_t u_r + \frac{n-1}{2}u u_t \\
& \qquad + \int_{K(S,T)} \frac{(n-1)^2}{2}\frac{\sin^2 u}{r^2} - \frac{(n-1)^2}{4}\,\frac{u\sin(2u)}{r^2} \\
& \qquad + \frac{1}{\sqrt{2}}\int_{C(S,T)} t(e^+ +m) + r (e^{-}+m) + \frac{n-1}{2}u (u_t + u_r)
\end{aligned}
\end{equation}

Our goal will be to show that
\begin{equation}
\lim_{T\to 0^+} \frac{1}{T}  \int_{\Sigma_T} Te^+ + r u_t u_r  = 0.
\label{redan}
\end{equation}
This implies \eqref{energyan}, which, based on Lemma \ref{biglem}, provides us also with the continuity of $u$. To achieve \eqref{redan}, we will prove that we can take $S=0$ in \eqref{energyint}, divide the resulting identity by $T$, and finally take the limit as $T\to 0$. 

\subsection{The surface integral} 
For $(t,r)\in C(S,T)$,  we notice that 
\begin{equation}
0 \leq  t(e^+ +m) + r (e^{-}+m) = t(u_t+u_r)^2 \lesssim t (e^+ +m)
\label{inf}
\end{equation}
Based on this, we prove:
\begin{lem}
\label{partiallem1}
If $u$ is a local smooth solution to \eqref{GWM} of finite energy, then 
\begin{equation}
\label{surfint}
 \int_{C(S,T)} t(e^+ +m) + r (e^{-}+m) + \frac{n-1}{2}\left|u (u_t + u_r)\right| \lesssim T F(T) + T^{\frac{n}{2}} F(T)^{\frac{1}{2}}.\end{equation}
\end{lem}
\begin{proof}
From \eqref{inf} and the definitions of $F(S,T)$, respectively $F(T)$, it is easy to see that
\[
0 \leq \int_{C(S,T)} t(e^+ +m) + r (e^{-}+m) \lesssim T F(S,T) \leq T F(T).
\]

For the last integrand in \eqref{surfint}, we can apply Cauchy-Schwarz inequality combined with the boundedness of $u$ to obtain
\[
\int_{C(S,T)}\left|u (u_t + u_r)\right| \lesssim \left(\int_{C(S,T)} u^2\right)^{\frac{1}{2}} \left(\int_{C(S,T)} (u_t+u_r)^2\right)^{\frac{1}{2}} \lesssim  \left(T^n - S^n\right)^{\frac{1}{2}}  F(S,T)^{\frac{1}{2}},\]
which provides the desired conclusion.
\end{proof}
 
We can immediately infer, using also the decay of flux \eqref{decayflux},
\begin{equation}
\lim_{T \to 0} \frac{1}{T} \left[\lim_{S \to 0} \int_{C(S,T)} t(e^+ +m) + r (e^{-}+m) + \frac{n-1}{2}u (u_t + u_r)\right] = 0,
\label{fxi}
\end{equation}
which finishes this analysis.

\subsection{The volume integrals}
First, for the volume integral on the left-hand side of \eqref{energyint}, we notice that the hypotheses on $g$ and $\alpha$ guarantees that both integrands are positive and so
\[\aligned
& \lim_{S \to 0}  \int_{K(S,T)} \left(\frac{n-1}{2}\right)^2 \frac{u f^{\prime}(u)f(u)}{r^{\alpha}} + (\alpha - (n+1))\frac{f^2(u)}{2r^{\alpha}}\\ 
& \qquad \qquad =  \int_{K(0,T)} \left(\frac{n-1}{2}\right)^2\frac{u f^{\prime}(u)f(u)}{r^{\alpha}} + (\alpha - (n+1))\frac{f^2(u)}{2r^{\alpha}}
\endaligned
\]
Next, for the right-hand side integral, we can prove:

\begin{lem}
\label{partiallem2}
If $u$ is a local smooth solution to \eqref{GWM} of finite energy and $T$ is sufficiently small, then 
\begin{equation}
\label{volint}
\int_{K(S,T)} \frac{\sin^2 u+|u\sin(2u)|}{r^2}  \lesssim 
\begin{cases}
T^{\frac{\alpha}{2}} & n=2,\\
T^{n-1} & n \geq 3.
\end{cases}
\end{equation}
\end{lem}
\begin{proof}
The easier of the two cases is when $n \geq 3$, for which the boundedness of $u$ implies 
\[
\begin{aligned}
\int_{K(S,T)} \frac{\sin^2 u+|u\sin(2u)|}{r^2} & = A(\mathbb{S}^{n-1})\int_S^T\int_0^t \frac{\sin^2 u+|u\sin(2u)|}{r^2} \ r^{n-1} drdt\\
&\lesssim \int_S^T\int_0^t\ r^{n-3}  drdt \lesssim  T^{n-1}
\end{aligned}
\]

The case $n=2$ is more subtle, as the continuity of $u$, obtained previously in Section 3.2, plays an important role.  First, as 
\[
\lim_{(t,r)\to (0,0)} u(t,r) =0,
\] 
$g(0) = 0$, and $f^{\prime}(0) \neq 0$, it follows that, for $(t,r)\in K(S,T)$ and $T$ sufficiently small,
\[
u^2 \sim u\sin (2u) \sim \sin^2 u \quad \text{and} \quad f(u) \sim u.\]
This allows us to deduce, based on the Cauchy-Schwarz inequality,
\[
u^2(t,r)  \lesssim \int_0^r |f(u(t,s))| |u_r(t,s)|\ ds   \lesssim  r^{\frac{\alpha-2}{2}} E(t).
\]

For $n=2$, $\alpha \geq 3$, which implies $\frac{\alpha-4}{2} \geq -\frac{1}{2}$. We can then estimate directly
\[
\int_{K(S,T)} \frac{\sin^2 u+|u\sin(2u)|}{r^2} \lesssim  \int_S^T\int_0^t \frac{u^2}{r^2} \ r drdt \lesssim \int_S^T\int_0^t\ r^{\frac{\alpha-4}{2}}  drdt \lesssim  T^{\frac{\alpha}{2}}.
\]

\end{proof}

\subsection{The time-slice integrals}
We make first the observation that
\[
S e^+ + r u_t u_r = (S-r)\frac{u_t^2+u_r^2}{2} + r \frac{(u_t+u_r)^2}{2} + S\left( \frac{n-1}{2}\, \frac{\sin^2 u(t,r)}{r^2} + \frac{f^2(u(t,r))}{2r^\alpha}\right)
\]
which tells us that the first two integrands combine to yield a non-negative quantity. Secondly, we show:

\begin{lem}
\label{partiallem3}
If $u$ is a local smooth solution to \eqref{GWM} of finite energy and $S$ is sufficiently small, then 
\begin{equation}
\int_{\Sigma_S} S e^+ + r \left|u_t u_r\right|   \lesssim S E(S),
\end{equation}
and
\begin{equation}
\int_{\Sigma_S} \left|u u_t\right| \lesssim
\begin{cases}
S^{\frac{\alpha+2}{4}} E(S), & n=2,\\
S^{\frac{n}{2}} E(S)^{\frac{1}{2}}, & n\geq 3.
\end{cases}
\label{uut}
\end{equation}
\end{lem}
\begin{proof}
For the first integrand, 
\[
\int_{\Sigma_S} S e^+  = S  E(S),\]
while the second one can be estimated using the Cauchy-Schwarz inequality,
\[
\int_{\Sigma_S} r \left|u_t u_r \right|  \leq S \left( \int_{\Sigma_S}u_t^2 \right)^{\frac{1}{2}} \left( \int_{\Sigma_S}u_r^2 \right)^{\frac{1}{2}} \leq S  E(S).
\]

Finally, for $n\geq 3$, the boundedness of $u$ implies
\[
\begin{aligned}
 \int_{\Sigma_S} \left|u u_t \right|  \lesssim \left( \int_{\Sigma_S}u^2 \right)^{\frac{1}{2}} \left( \int_{\Sigma_S}u_t^2 \right)^{\frac{1}{2}}\
& \lesssim \left( \int_0^S u^2(S,r)\ r^{n-1}dr \right)^{\frac{1}{2}} E^{\frac{1}{2}}(S)\\
& \lesssim S^{\frac{n}{2}} E(S)^{\frac{1}{2}}.
\end{aligned}
\]
The argument for the $n=2$ case is similar to the one above, but one has to use also, as for the volume integrals, the rate of decay of $u$.
\end{proof}

\subsection{Conclusion for the proof of Theorem \ref{bigthm1}}
We now have all the ingredients necessary for finishing the main argument. First, using the results contained in Lemmas \ref{partiallem1}-\ref{partiallem3}, we can take the limit as $S \to 0$ in \eqref{energyint} to obtain
\begin{equation}
\label{partialestimate}
\begin{aligned}
& \int_{K(0,T)} \frac{(n-1)^2}{4}\frac{u\, f^{\prime}(u)f(u)}{r^{\alpha}} + (\alpha - (n+1))\frac{f^2(u)}{2r^{\alpha}} \qquad \qquad \qquad \\
& + \int_{\Sigma_T} Te^+ + r u_t u_r  = \int_{K(0,T)} \frac{(n-1)^2}{2}\frac{\sin^2 u}{r^2} - \frac{(n-1)^2}{4}\frac{u\sin(2u)}{r^2} \\
& \qquad + \frac{1}{\sqrt{2}}\int_{C(0,T)} t(e^+ +m) + r (e^{-}+m) + \frac{n-1}{2}\,u \,(u_t + u_r) \\
& \qquad - \int_{\Sigma_T} \frac{n-1}{2}\, u \,u_t
\end{aligned}
\end{equation}
Next, we divide by $T$ and take the limit as $T \to 0$. Taking advantage of \eqref{fxi}, \eqref{volint}, and \eqref{uut}, we deduce
\begin{equation}
\aligned
& \lim_{T\to 0} \frac{1}{T}\bigg[ \int_{K(0,T)} \frac{(n-1)^2}{4}\frac{u f^{\prime}(u)f(u)}{r^{\alpha}} + (\alpha - (n+1))\frac{f^2(u)}{2r^{\alpha}}\\
& \qquad \qquad \qquad \qquad +  \int_{\Sigma_T} Te^+ + r u_t u_r \bigg]\,=\,0.
\endaligned
\label{fin}
\end{equation}

Using previously made remarks concerning the positivity of integrands in \eqref{fin}, we conclude that \eqref{redan} holds, thus finishing the proof.

\section{Acknowledgements}
The author was supported by National Science Foundation Career Grant DMS-0747656.

\bibliographystyle{amsplain}
\bibliography{database}

\providecommand{\bysame}{\leavevmode\hbox to3em{\hrulefill}\thinspace}
\providecommand{\MR}{\relax\ifhmode\unskip\space\fi MR }
\providecommand{\MRhref}[2]{%
  \href{http://www.ams.org/mathscinet-getitem?mr=#1}{#2}
}
\providecommand{\href}[2]{#2}
\begin{thebibliography}{1}

\bibitem{Adkins1984}
Gregory Adkins and Chiara Nappi, \emph{Stabilization of chiral solitons via
  vector mesons}, Phys. Lett. B \textbf{137} (1984), no.~3-4, 251--256.

\bibitem{Derrick1964}
G.~H. Derrick, \emph{Comments on nonlinear wave equations as models for
  elementary particles}, J. Mathematical Phys. \textbf{5} (1964), 1252--1254.
  \MR{0174304 (30 \#4510)}

\bibitem{Geba2010}
Dan-Andrei Geba and Sarada~G. Rajeev, \emph{A continuity argument for a
  semilinear {S}kyrme model}, Electron. J. Differential Equations (2010), No.
  86, 9. \MR{2680289 (2011f:83068)}

\bibitem{Morawetz1975}
Cathleen Morawetz, \emph{Notes on time decay and scattering for some hyperbolic
  problems}, Regional Conference Series in Applied Mathematics, vol.~19,
  Society for Industrial and Applied Mathematics, Philadelphia, 1975.

\bibitem{Shatah1992}
J.~Shatah and A.~Tahvildar-Zadeh, \emph{Regularity of harmonic maps from the
  {M}inkowski space into rotationally symmetric manifolds}, Comm. Pure Appl.
  Math. \textbf{45} (1992), no.~8, 947--971.

\end{thebibliography}

\end{document}